\documentclass{gtart_a}
\pdfoutput=1

\title{Generic uniqueness of least area planes in hyperbolic space}

\author{Baris Coskunuzer}
\givenname{Baris}
\surname{Coskunuzer}
\address{Department of Mathematics\\
Yale University\\\newline
New Haven CT 06520\\
USA}
\email{baris.coskunuzer@yale.edu}
\urladdr{}

\volumenumber{10}
\issuenumber{}
\publicationyear{2006}
\papernumber{12}
\lognumber{0484}
\startpage{401}
\endpage{412}

\doi{}
\MR{}
\Zbl{}

\arxivreference{math.GT/0408066}
\arxivpassword{qdacf}   


\received{5 August 2004}
\revised{23 February 2006}
\accepted{27 March 2006}
\publishedonline{27 April 2006}
\published{27 April 2006}
\proposed{Gang Tian}
\seconded{Tobias Colding, David Gabai}
\corresponding{}
\editor{}
\version{}

\AtBeginDocument{\let\bar\wbar\let\hat\what}

\newtheorem{thm}{Theorem}[section]
\newtheorem*{mainthm}{Main theorem} 
\newtheorem{lem}{Lemma}[section]
\makeatletter
\let\c@lem\c@thm
\makeatother

\theoremstyle{definition}
\newtheorem{defn}{Definition}[section]

\theoremstyle{remark}
\newtheorem{rmk}{Remark}[section]

\newcommand{\Si}{S^2_{\infty}(\mathbb{H}^3)}
\newcommand{\PI}{\partial_{\infty}}
\newcommand{\BH}{\mathbb{H}}
\newcommand{\BR}{\mathbb{R}}

\makeop{sup}
\makeop{inf}

\begin{document}

\begin{asciiabstract}
We study the number of solutions of the asymptotic Plateau problem in
H^3.  By using the analytical results in our previous paper, and some
topological arguments, we show that there exists an open dense subset
of C^3 Jordan curves in S^2_{infinity}(H^3) such that any curve in this
set bounds a unique least area plane in H^3.
\end{asciiabstract}

\begin{htmlabstract}
We study the number of solutions of the asymptotic Plateau
problem in <b>H</b><sup>3</sup>. By using the analytical results
in our previous paper, and some topological arguments, we show that
there exists an open dense subset of C<sup>3</sup> Jordan curves in
S<sup>2</sup><sub>&infin;</sub>(<b>H</b><sup>3</sup>) such that any curve
in this set bounds a unique least area plane in <b>H</b><sup>3</sup>.
\end{htmlabstract}

\begin{abstract}
We study the number of solutions of the asymptotic Plateau problem in
$\mathbb{H}^3$. By using the analytical results in our previous paper,
and some topological arguments, we show that there exists an open dense
subset of $C^3$ Jordan curves in $S^2_{\infty}(\mathbb{H}^3)$ such that
any curve in this set bounds a unique least area plane in $\mathbb{H}^3$.
\end{abstract}

\maketitle

\section{Introduction}
\label{sec1}

In this paper, we are interested in the number of solutions to the
asymptotic Plateau problem. The asymptotic Plateau problem was solved
by Michael Anderson in his seminal paper \cite{A1}.  Furthermore in
\cite{A2} he used geometric measure theory methods to prove that for
any Jordan curve in $\Si$, there exists an embedded least area plane
in $\BH^3$ spanning that curve.  Then, Hardt and Lin studied the
asymptotic regularity of the least area surfaces in hyperbolic space
in \cite{HL,Li1}. Also, Gabai \cite{Ga} solved the asymptotic Plateau
problem by using topological arguments developed by Hass and Scott
\cite{HS}. Lang \cite{L} generalized Anderson's methods to solve the
problem in Hadamard Gromov hyperbolic spaces. In \cite{Co2}, the
author solved the problem for Gromov hyperbolic spaces with cocompact
metric by generalizing Gabai's techniques.

On the number of the least area planes spanning a given asymptotic curve,
there are a few results so far.  Anderson \cite{A1} showed that if the
curve in $\Si$ bounds a convex domain in $\Si$, then there exists a unique
least area plane in $\BH^3$ spanning that curve. Then Lin \cite{Li1}
generalized this result to the curves bounding a star shaped domain
in $\Si$. Recently, in \cite{Co1}, the author showed that for a generic
$C^{3,\mu}$ Jordan curve in $\Si$, there exist finitely many least area
planes in $\BH^3$ asymptotic to that curve. In this paper, we improve
the generic finiteness result of the previous paper \cite{Co1} to a
generic uniqueness result for a more general class of curves.

Here is an outline of the paper.  In \cite{Co1}, we showed that the
boundary restriction map $\pi$ from the space of minimal maps from
$D^2$ to $\BH^3$ with $C^3$ asymptotic data to the space of the $C^3$
immersions of $S^1$ into $\Si$ is Fredholm of index $0$. We consider
the boundary restriction map $\pi$ from minimal maps to their boundary
parametrizations. Since it is Fredholm of index $0$, its derivative is
an isomorphism for generic curves.

Fix a generic curve $\Gamma$ in $\Si$. By using the inverse function
theorem, we find a neighborhood $U_\Sigma$ of a least area plane $\Sigma$
in $\pi^{-1}(\Gamma)$, mapping homeomorphically into a neighborhood
$V_\Gamma$ of $\Gamma$.  By taking a path $\alpha$ in $V_\Gamma$, and
by considering the corresponding path $\pi^{-1}(\alpha)$ in $U_\Sigma$,
we get a continuous family of minimal planes with disjoint asymptotic
boundaries around $\Sigma$. Then, we show that this continuous family
of minimal planes is indeed a foliation by least area planes of a
neighborhood of $\Sigma$. This implies the uniqueness of the least area
plane in $\BH^3$ spanning $\Gamma$. Then we show the same result for
any curve in a neighborhood of a generic curve, and we get an open dense
subset of the $C^3$ Jordan curves in $\Si$ with the uniqueness result.

\begin{mainthm}
Let $A$ be the space of $C^3$ simple closed curves in $\Si$. Then there
exists an open dense subset $A'$ in $A$ such that for any $\Gamma \in A'$,
there exists a unique least area plane $\Sigma$ with $\PI \Sigma=\Gamma$.
\end{mainthm}

The organization of the paper is as follows. In the next section, we
will give some preliminary results which will be used in the following
sections. In \fullref{sec3}, we will construct a foliated neighborhood
of a least area plane spanning a generic curve. By using this, we will
show that there is a unique least area plane spanning that curve. In
\fullref{sec4}, we will show that the curves with this property is
indeed open dense subset in the ambient space. Finally in
\fullref{sec5} we will give some concluding remarks.

\subsection{Acknowledgements}
I want to thank David Gabai for very useful conversations, and helpful
remarks. I want to thank the referee for very valuable remarks.

\section{Preliminaries}
\label{sec2}

In this section, we will overview the basic results which we use in the
following sections.

\begin{defn}
A {\em minimal plane} is a plane such that the mean curvature is $0$ at
every point. A {\em least area plane} is a plane such that any disk in
the plane has least area among the disks with same boundary. A least area
plane is minimal, but the converse is not true in general.
\end{defn}

\begin{defn}
A linear operator between two Banach spaces is called {\em a Fredholm
operator} if the dimension of the kernel and codimension of the image are
finite. The difference between these dimensions is called the {\em
Fredholm index} of the operator. A map between two Banach manifolds is
{\em a Fredholm map} if the derivative is a Fredholm operator at every
point. A {\em regular value} is a point in the image such that for any
point in the preimage (may be empty), the derivative is surjective.
\end{defn}

The classical theorem of the subject is the Sard--Smale theorem from
\cite{Sm}.

\begin{thm}[Sard--Smale]
Let $f\co X\rightarrow Y$ be a Fredholm map. Then the regular values of $f$
are almost all of Y, that is, except a set of
the first category.
\end{thm}

We will use the following spaces in the remaining part of the paper.
\begin{eqnarray*}
A &=& \bigl\{\alpha \in C^3(S^1,S^2_\infty) \ | \ \alpha
  \text{ embedding }\bigr\} \\
D &=& \bigl\{u \in C^3(S^1,S^1) \ | \  u \text{ diffeomorphism and }
u(e^{2k\pi i/3}) = e^{2k\pi i/3}, k=1,2,3 \bigr\} \\
M &= & \bigl\{f\co D^2\rightarrow \BH^3 \ | \ f(D^2) \text{ minimal and }
  f|_{\partial D^2} \in A \bigr\}
\end{eqnarray*}
Now, we will summarize the results of \cite{Co1} which are essential
for this paper. In \cite{Co1}, we studied the asymptotic Plateau problem
by using global analysis methods. We first looked at the space $M$. We
identified the minimal maps in $M$ with the conformal harmonic maps,
and considered $M$ as a subspace of the space of harmonic maps. Then by
using Li and Tam's results \cite{LT1,LT2}, we identified these harmonic
maps with their boundary parametrizations, corresponding to the space
$A \times D$.

In other words, we realize $M$, the space of minimal planes, as a subspace
of the product bundle $A \times D$, by identifying $(\alpha,u)$ in $A$
with its unique harmonic extension $\widetilde{\alpha{\circ}u}$. Let
$\alpha \in A$, and consider the fiber $(\{\alpha\}\times D )\cap M$. This
set contains all minimal planes, $\Sigma_i=\widetilde{\alpha{\circ}
u_i}(D^2)$, with asymptotic boundary $\Gamma= \alpha(S^1)$, that is, $\PI
\Sigma_i=\Gamma$.  Here, the projection map $\pi$ from $A\times D$
into $A$, sends the elements of $M$ to their boundary parametrizations.
By using Tomi and Tromba's techniques in \cite{TT}, we proved that $\pi
{|_M}$ is Fredholm of index $0$.

\begin{thm}
\label{thm22}
$M=\{(\alpha,u)\in A \times D \ |\  \widetilde{\alpha{\circ}u}(D^2)
\mbox{ is minimal}\, \}$ is a submanifold of the product bundle $A \times D$; and
the bundle projection map when restricted to $M$, $\pi |_M\co M\rightarrow
A$, is Fredholm of index 0.  \end{thm}

Then by using the boundary regularity results of Hardt and Lin
\cite{HL} to get compactness, we showed that for a generic
curve $\Gamma$ in $A$, there are finitely many least area planes in
$\pi^{-1}(\Gamma)$. In other words, we proved that for a generic curve,
there exist finitely many least area planes spanning that curve. Indeed,
this is true in more generality. By using the more general regularity
results of Tonegawa \cite{To} instead of those of \cite{HL}, the same
technique shows that for a generic curve, there exist finitely many {\it
minimal planes} spanning that curve.

\begin{thm}
Let $\Gamma \subset \Si$ be a generic curve in $A$. Then there are finitely many minimal planes $\Sigma$ in $\BH^3$
with $\PI \Sigma =\Gamma$.
\end{thm}

\begin{proof}
In \cite[Theorem~1.6]{To}, Tonegawa shows that if the asymptotic curve
$\Gamma\subset \Si$ is $C^{3,\mu}$, and the minimal plane $\Sigma\subset
\BH^3$ ($H=0$) asymptotic to $\Gamma$, then $\Sigma \cup \Gamma$ is
also $C^{3,\mu}$ regular up to the boundary. In \cite[Lemma~4.5]{Co1},
we used the same result for least area planes from \cite{Li1}. If we
replace this boundary regularity result for least area planes with
the more general result for minimal planes, we see that 
\cite[Lemma~4.5]{Co1} is true for minimal planes, too. Then this implies by 
\cite[Theorem~4.7]{Co1}, there are finitely many minimal planes for a generic
asymptotic curve.
\end{proof}

\begin{rmk}
With the above notation, this theorem implies that for a generic curve $\alpha \in A$, $\pi ^{-1}(\alpha) \subset M$ is
a finite set, that is, generic curves have finitely many preimages in $M$
under the projection map $\pi \co M\rightarrow A$.
\end{rmk}

\begin{rmk}
In \cite{Co1}, we took $A$ as the space of $C^3$ immersions instead of $C^3$ embeddings. In this paper, we are only interested
in $C^3$ embeddings. Since $A$ is an open submanifold of the space of immersions, we can replace it with the
corresponding space in \cite{Co1}, and the same proofs will go through.
\end{rmk}

The next theorem we use is the classical inverse function theorem for
Banach manifolds (see Lang \cite{La}).

\begin{thm}[Inverse Function Theorem] Let $M$ and $N$ be Banach manifolds,
and let $f\co M \rightarrow N$ be a $C^p$ map. Let $x_0
\in M$ and $df$ be an isomorphism at $x_0$. Then $f$ is local $C^p$
diffeomorphism, that is, there exists an open neighborhood
of $U\subset M$ of $x_0$ and an open neighborhood $V \subset N$ of
$f(x_0)$ such that $f|_U\co U\rightarrow V$ is $C^p$
diffeomorphism.
\end{thm}

Now, we can establish the starting point of the paper.

\begin{thm}
\label{thm25}
Let $\Gamma \in A$ be a generic curve as above. Then for any $\Sigma\in \pi^{-1}(\Gamma)$, there exist neighborhoods
$U_\Sigma\subset M$, and $V_\Gamma\subset A$ such that $\pi
|_{U_\Sigma}\co U_\Sigma \rightarrow V_\Gamma$ is a
homeomorphism.
\end{thm}

\begin{proof}
By \fullref{thm22}, the map $\pi |_M\co M \rightarrow A$ is Fredholm of
index 0. By the Sard--Smale Theorem, the regular values are generic for
$\pi$. By the work of Anderson \cite{A2}, we also know that $\pi |_M$
is surjective. Let $\Gamma\in A$ be a regular value, and $\Sigma \in
\pi^{-1}(\Gamma)\subset M$. Since $\Gamma$ is a regular value, $D\pi
(\Sigma)\co T_\Sigma M \rightarrow T_\Gamma A$ is surjective. Moreover,
we know that $\pi$ is Fredholm of index $0$. This implies $D\pi$ is an
isomorphism at the point $\Sigma\in M$. By the Inverse Function Theorem,
there exists a neighborhood $U_\Sigma$ of $\Sigma$ in $M$, which $\pi$
maps homeomorphically onto a neighborhood $V_\Gamma$ of $\Gamma$ in $A$.
\end{proof}

\section{Foliated neighborhoods of least area planes}
\label{sec3}

In this section, our aim is to construct a foliated neighborhood for any least area plane spanning a given generic
curve in $A$. Moreover, we will show that the leaves of this foliation are embedded least area planes whose asymptotic
boundaries are disjoint from each other. By using this, we will show the uniqueness of the least area plane spanning
the generic curve.

We will abuse the notation by using interchangeably the map $\Gamma\co S^1\rightarrow \Si$ with its image $\Gamma(S^1)$.
Similarly the same is true for $\Sigma \co D^2\rightarrow \BH^3$ and its image $\Sigma(D^2)$.

Let $\Gamma_0\in A$ be a generic curve, and let $\Sigma_0\in \pi ^{-1}(\Gamma_0)\subset M$ be a least area plane whose
existence is guaranteed by Anderson's result in \cite{A2}. Then by
\fullref{thm25}, there is a neighborhood of $\Sigma_0\in
U\subset M$ that is homeomorphic to the neighborhood $\Gamma_0\in V \subset A$.

Let $\Gamma\co [-\epsilon,\epsilon]\rightarrow V$ be a path such that $\Gamma(0)=\Gamma_0$ and for any $t, t'
\in[-\epsilon,\epsilon]$, $\Gamma_t\cap\Gamma_{t'}=\emptyset$. In other words, $\{\Gamma_t\}$ foliates a neighborhood
of $\Gamma_0$ in $\Si$. Let $\Sigma_t\in U$ be the preimage of $\Gamma_t$ under the homeomorphism.

\begin{lem}
\label{lem31}
$\{\Sigma_t\}$ is a foliation of a neighborhood of $\Sigma_0$ in $\BH^3$ by embedded least area planes.
\end{lem}

\begin{proof} We claim that we can assume that $\Sigma_{\pm\epsilon}$ are least area planes. By Remark
2.1, we know that $\Gamma_0$ has finitely many preimages, say $\{S_1,S_2,
... , S_k\}\subset M$. Then by \fullref{thm25},
for each $S_i$, we have a neighborhood $U_i\subset M$ such that $\pi
|_{U_i}\co U_i\rightarrow V_i$ is a homeomorphism
where $V_i\subset A$ is a neighborhood of $\Gamma_0$. Now, by taking the intersection, we can assume $V_i=V_j$ (say
$V$) and $\pi ^{-1}(V)=\{U_1, U_2, ..., U_k\}$. We also assume that $\{\Gamma_t\}\subset V$. By Anderson's theorem, we
know that there exists a least area plane $\Sigma_t$ for any $\Gamma_t$ with $\PI\Sigma_t=\Gamma_t$. This means for any
$t\in(-\epsilon, \epsilon)$, there exists a least area plane $\Sigma_t\in \pi ^{-1}(\Gamma_t)$, and it must belong to
$U_i$ for some $i$. Since there are infinitely many $\Sigma_t$, we can
find a $U_i$ such that there are three least area
planes $\Sigma_{t_0}, \Sigma_{t_1},\Sigma_{t_2}$ in the preimage of
$\{\Gamma_t\}$ restricted to $U_i$, that is,
$\Sigma_{t_j}\in U_i\cap\pi ^{-1}(\{\Gamma_t\})$. Now let $t_0<0 \leq t_1<t_2$, and we can reparametrize the path
$\Gamma$ so that $\Gamma_{-\epsilon}=\Gamma_{t_0}$, $\Gamma_0=\Gamma_{t_1}$, $\Gamma_{+\epsilon}=\Gamma_{t_2}$. From
now on, we will assume that $\Sigma_{\pm\epsilon}$ are also least area planes.

We will prove the lemma in three steps.

\textbf{Step 1}\qua For any $s\in [-\epsilon, \epsilon]$, $\Sigma_s$ is an embedded plane.

Since $\Sigma_0$ is a least area plane, by \cite[Theorem~4.5]{Co1},
$\Sigma_0$ is an embedded plane. Now, $\{\Sigma_t\}$ is a continuous
family of minimal planes. We cannot apply the same theorem to these
planes, since the theorem is true for least area planes, while our planes
are only minimal.

Let $s_0=\inf\{s\in(0,\epsilon] \ | \ \Sigma_s \mbox{ is not embedded}\}$. However, since $\{\Sigma_t\}$ is continuous
family of planes, and this can only happen when $\Sigma_{s_0}$ has tangential self intersection (locally lying on one
side), this contradicts with the maximum principle for minimal surfaces. So for all $s\in [0,\epsilon]$, $\Sigma_s$ is
embedded. Similarly, this is true for $s\in [-\epsilon, 0]$, and the result follows.

\textbf{Step 2}\qua $\{\Sigma_t\}$ is a foliation, that is, for any $t, t' \in[\epsilon,\epsilon]$,
$\Sigma_t\cap\Sigma_{t'}=\emptyset$.

Assume on the contrary that there exist $t_1 < t_2$ such that $\Sigma_{t_1}\cap\Sigma_{t_2}\neq \emptyset$. First,
since the asymptotic boundaries $\Gamma_{t_1}$ and $\Gamma_{t_2}$ are disjoint, the intersection cannot contain an
infinite line. So the intersection must be a collection of closed curves. We will show that in this situation, there
must be a tangential intersection between two planes, and this will contradict to the maximum principle for minimal
surfaces.

If $\Sigma_{t_2}$ does not intersect all the minimal planes $\Sigma_s$ for
$s\in [-\epsilon, t_2]$, let
$$s_0=\sup\{s\in [-\epsilon, t_2] \ | \ \Sigma_{t_2}\cap\Sigma_s =
\emptyset\}.$$
Then, since $\{\Sigma_t\}$ is a continuous family of
minimal planes, it is clear that $\Sigma_{t_2}$ must intersect $\Sigma_{s_0}$ tangentially, and lie in one side of
$\Sigma_{s_0}$. However, this contradicts with the maximum principle for minimal surfaces.

So, let's assume $\Sigma_{t_2}$ intersects all minimal planes $\Sigma_s$ for $s\in [-\epsilon, t_2]$. Let
$$s_0=\sup\{s\in(-\epsilon,t_2] \ | \ \Sigma_{-\epsilon} \cap \Sigma_s
=\emptyset\}.$$
If there exists an $s_0\in
(-\epsilon, t_2]$, then again since $\{\Sigma_t\}$ is a continuous family of minimal planes, $\Sigma_{s_0}$ must
intersect $\Sigma_{-\epsilon}$ tangentially, and lie in one side of $\Sigma_{s_0}$. However, this contradicts with the
maximum principle for minimal surfaces.

If there is no such $s_0\in (-\epsilon, t_2]$, then this means $\Sigma_0\cap \Sigma_{-\epsilon} \neq \emptyset$. Since
$\PI \Sigma_0 =\Gamma_0$ is disjoint from $\PI \Sigma_{-\epsilon} =\Gamma_{-\epsilon}$ in $\Si$, the intersection
cannot contain an infinite line. So, we can find a simple closed loop $\gamma$ in the intersection $\Sigma_0\cap
\Sigma_{-\epsilon}$. But, since $\Sigma_0$ and $\Sigma_{-\epsilon}$ are both least area planes, this is a contradiction
by the Meeks--Yau exchange roundoff trick.

\textbf{Step 3}\qua For any $s\in [-\epsilon, \epsilon]$, $\Sigma_s$ is a least area plane.

Fix $\Sigma_s$ for $s\in (-\epsilon, \epsilon)$. We will show that any subdisk of $\Sigma_s$ is a least area disk. Let
$[\Sigma_{-\epsilon}, \Sigma_{\epsilon}]$ be the region bounded by embedded planes $\Sigma_{-\epsilon}$ and
$\Sigma_{\epsilon}$ in $\BH^3$. By above results, $\Sigma_s\subset [\Sigma_{-\epsilon}, \Sigma_{\epsilon}]$. Let
$\gamma\subset \Sigma_s$ be an extreme (in the boundary of its convex hull) simple closed curve. By Morrey's theorem
\cite{Mo}, there exists a least area disk $D$ with $\partial D=\gamma$. By
the work of Anderson \cite{A1}, the least area disks in hyperbolic space have
convex hull property, so $D$ is in its convex hull. Since we chose $\gamma$ to be extreme, by Meeks and Yau's theorem
\cite{MY}, $D$ is a least area embedded disk in $\BH^3$ spanning $\gamma$.

We claim that $D\subset [\Sigma_{-\epsilon}, \Sigma_{\epsilon}]$. Otherwise, $D$ must intersect one or both of the
planes. Since the boundary of $D$, is disjoint from $\Sigma_{\pm\epsilon}$, then there must be a simple closed curve in
the intersection. But all of them are least area, and this is a contradiction by exchange roundoff trick of Meeks and
Yau.

Now, we want to show that $D\subset \Sigma_s$. If $D$ is not contained in $\Sigma_s$, it must intersect other leaves
nontrivially. Since $D$ is completely in the domain $[\Sigma_{-\epsilon}, \Sigma_{\epsilon}]$, then $\{\Sigma_t\}\cap
D$ induce a singular 1--dimensional foliation $F$ on $D$. The singularities of the foliation are isolated, as
$\{\Sigma_t\}$ are minimal planes. Since Euler characteristic of the disk
is 1, by the Poincar\'e--Hopf index formula there
must be a positive index singularity implying tangential (lying on one side) intersection of $D$ with some leave
$\Sigma_s$. However, this contradicts with the maximum principle for minimal surfaces. So, we show that the subdisks of
$\Sigma_s$ with extreme boundary are least area disks. But since for any subdisk of $\Sigma_s$, we can find a bigger
disk containing it with an extreme boundary, it is true for any subdisk of $\Sigma_s$. This shows that for any $s\in
(-\epsilon,+\epsilon)$, $\Sigma_s$ is a least area plane.
\end{proof}

\begin{lem}
\label{lem32}
Let $\Gamma$ be a generic curve as described above. Then, there exists a unique least area plane $\Sigma$ with $\PI
\Sigma = \Gamma$.
\end{lem}

\begin{proof} We will use the notation of the proof of the previous lemma. In the proof of the previous lemma, we
constructed a foliation by least area planes $\{\Sigma_t\}$ with asymptotic boundaries $\{\Gamma_t\}$. After the
reparametrization in the beginning of the proof of the previous lemma, the generic curve $\Gamma$ becomes an interior
leaf of the foliation $\{\Gamma_t\}$. In other words, there is a $t'\in (-\epsilon,\epsilon)$ such that $\Gamma =
\Gamma_{t'}$. Since $t'$ is an interior point, $\{\Sigma_t\}$ is also foliation for a neighborhood of $\Sigma_{t'}$.

Let $\widehat{\Sigma}$ be another least area plane with asymptotic boundary $\Gamma_{t'}$. If $\widehat{\Sigma} \neq
\Sigma_{t'}$ then $\widehat{\Sigma}$ must intersect a leave in the foliated neighborhood of $\Sigma_{t'}$, say
$\Sigma_s$. But, since $\PI\Sigma_s = \Gamma_s$ is disjoint from $\PI \Sigma_{t'}=\Gamma_{t'}$, this implies the
intersection cannot have infinite lines, but closed loops. However, these
are least area planes, and by the Meeks--Yau
exchange roundoff trick, two least area planes cannot intersect in a closed loop. This is a contradiction.
\end{proof}

\begin{rmk} This lemma does not say that there exists a unique {\it minimal} plane spanning a given generic curve.
In the proof of the lemma, we essentially use the plane being least area.
\end{rmk}

Hence, we have proved the following theorem:

\begin{thm}
Let $\Gamma$ be a generic curve in $A$. Then there exists a unique least area plane $\Sigma$ in $\BH^3$ with
$\PI\Sigma=\Gamma$.
\end{thm}

So far we have proved the uniqueness of least area planes for a subset $\hat{A}\subset A$ which is except a first
category set. In the next section, we will show that this is true for a
more general class of curves, that is, an open
dense subset of $A$.

\section{Open dense set of curves}
\label{sec4}

In this section, we will show that any generic curve, which is a regular value of the Fredholm map $\pi$, has an open
neighborhood such that the uniqueness result holds for any curve in this neighborhood.

Let $\Gamma_0\in A$ be a generic curve, and let $\Sigma_0\in \pi ^{-1}(\Gamma_0)\subset M$ be the unique least area
plane spanning $\Gamma_0$. Let $ U\subset M$ be the neighborhood of $\Sigma_0$ homeomorphic to the neighborhood $V
\subset A$ of $\Gamma_0$ as in the previous section. We will show that $\Gamma_0$ has a smaller open neighborhood
$V'\subset V$ such that for any $\Gamma\in V'$, there exists a unique least area plane in $\BH^3$ with
$\PI\Sigma=\Gamma$.

First we will show that the curves disjoint from $\Gamma_0$ in the open neighborhood also bounds a unique least area
plane in $\BH^3$.

\begin{lem} Let $\beta\in V$ with $\beta\cap\Gamma_0= \emptyset$. Then
there exists a unique least area plane $\Sigma$
in $\BH^3$ with $\PI \Sigma = \beta$.
\end{lem}

\begin{proof}
Since $\beta\in V$ is disjoint from $\Gamma_0$, we can find a path
$\Gamma\co (-\epsilon,\epsilon)\rightarrow V$, such
that $\{\Gamma_t\}$ foliates a neighborhood of $\Gamma_0$ in $\Si$, and
$\beta$ is one of the leaves, that is,
$\beta=\Gamma_s$ for some $s\in(-\epsilon,\epsilon)$. Then the proofs of the previous section implies that
$\Sigma_\beta=\Sigma_s$ and $\{\Sigma_t\}$ also gives a foliation of a neighborhood of $\Sigma_\beta$ by least area
planes. Then proof of \fullref{lem32} implies that $\Sigma_\beta$ is the unique least area plane spanning $\beta$.
\end{proof}

Now, if we can show the same result for the curves in $V$ intersecting $\Gamma_0$, then we are done. Unfortunately, we
cannot do that, but we will bypass this by going to a smaller neighborhood.

\begin{lem} There exists a neighborhood $V'\subset V$ of $\Gamma_0$ such
that for any $\Gamma_0'\in V'$, there exists a
unique least area plane with asymptotic boundary $\Gamma_0'$.
\end{lem}

\begin{proof}
Let $\beta_1, \beta_2 \in V$ be two Jordan curves in $\Si$ disjoint from
$\Gamma_0$ and lying in opposite side of $\Gamma_0$ in $\Si$. Then,
$\beta_1$ and $\beta_2$ bounds an open annulus $B$ in $\Si$, where
$\Gamma_0$ is in the interior of $B$. Then, let $V' = \{\Gamma \in V
\ | \Gamma \subset B \}$ be an open neighborhood of $\Gamma_0$.

Now, fix $\Gamma_0'\in V'$. By the assumption on
$V'$, $\Gamma_0,\Gamma_0'$ both belong to the annulus
$B\subset \Si$.  Then, we can find two paths $\Gamma,\Gamma'\co
[-\epsilon,\epsilon]\rightarrow V$ where $\{\Gamma_t\}$, $\{\Gamma_t'\}$
foliate $B$ such that $\Gamma(\epsilon)={\Gamma'}(\epsilon)=\beta_1$
, $\Gamma({-\epsilon})={\Gamma'}({-\epsilon})=\beta_2$, and
$\Gamma(0)=\Gamma_0$, $\Gamma'(0)=\Gamma_0'$.

By \fullref{lem31}, we know that $\{\Gamma_t\}$ induces $\{\Sigma_t\}$ family of embedded least area planes with asymptotic
boundary $\{\Gamma_t\}$. Moreover, these least area planes are unique with the given asymptotic boundary, and the
leaves of the foliation in the neighborhood of $\Sigma_0$.

Now, consider the preimage of the path $\Gamma'$ under the
homeomorphism $\pi |_U\co U \rightarrow V$. This will give us a
path $\Sigma' \subset U\subset M$, which is a continuous family
of minimal planes, say $\{\Sigma_t'\}$. We claim that this is
also a family of embedded least area planes inducing a foliated
neighborhood of $\Sigma_0'$. By previous section, we can also
assume that $\beta_1$ and $\beta_2$ bound a unique least area plane,
say $\PI \Sigma_\epsilon = \beta_1$ and $\PI \Sigma_{-\epsilon} =
\beta_2$. This means $\Sigma_{\pm\epsilon}'=\Sigma_{\pm\epsilon}$. So,
the family $\{\Sigma_t'\}$ has embedded least area planes
$\Sigma_{\pm\epsilon}'$. Then a slight modification of the proof of
\fullref{lem31} shows that $\{\Sigma_t'\}$ is a family of embedded least
area planes inducing a foliation of a neighborhood of $\Sigma_0'$. By
\fullref{lem32}, $\Sigma_0'$ is the unique least area plane spanning
$\Gamma_0'$.
\end{proof}

Thus, we get the following theorem.

\begin{thm}
There exists an open dense subset $A'$ in $A$, where $$A= \{\alpha
\in C^3(S^1,S^2_\infty) \ | \ \alpha \text{ embedding}\}$$ with $C^3$
topology, such that for any $\Gamma\in A'$, there exists a unique least
area plane with asymptotic boundary $\Gamma$.
\end{thm}

\begin{proof}
The set of regular values of Fredholm map, say $\hat{A}$, is the whole
set except a set of the first category by the Sard--Smale theorem. So, the
regular curves are dense in $A$. By the above lemmas, for any regular curve
$\Gamma_0$, there exists an open neighborhood $V_{\Gamma_0}'\subset A$
which the uniqueness result holds. So, $A'=\bigcup_{\Gamma \in\hat{A}}
V_\Gamma'$ is an open dense subset with the desired properties.
\end{proof}

\section{Final remarks}
\label{sec5}

\subsection{The technique}

The technique we used in this paper to show uniqueness of the least
area planes is similar to the techniques of Anderson and Lin for their
proofs for convex domains \cite{A1} and star shaped domains \cite{Li1}
in $\Si$. The idea is simple. Fix the least area plane. Construct a
neighborhood of the least area plane foliated by least area planes with
disjoint asymptotic boundaries. Then any other least area plane with
same asymptotic boundary with the original one must intersect the other
least area planes. But this intersection cannot contain an infinite line
as the asymptotic boundaries are disjoint. So the intersection must be
a collection of closed curves. But the Meeks--Yau exchange roundoff trick
says that two least area planes cannot intersect in a closed loop,
and gives the desired contradiction.

On the other hand, hyperbolic space and its asymptotic boundary are not
very essential for this technique. One can employ the same method for any
convex domain, and least area disks whose boundaries are in the boundary
of this convex domain, once one has a continuous family of minimal disks.

\subsection{Extremal curves in $\BR^3$}

As we mentioned above, we can use this technique for extremal curves
in convex domains in $\BR^3$. Let $N\subset\BR^3$ be a smooth convex
domain. Then by changing the space of curves from embeddings of $S^1$
into $\BR^3$ to the embeddings of $S^1$ into $\partial N$ in the paper by
Tomi and Tromba \cite{TT},
the whole proof would go through. So this will establish the necessary
analytical background to employ the techniques in this paper. Then one
can fix a regular extremal curve in $\partial N$, and get a foliated
neighborhood of the embedded least area disk spanning the fixed curve
in $N$. This proves again the uniqueness of the least area disk. So,
one can show that for a given smooth convex domain $N$ in $\BR^3$,
there exists an open dense subset of the Jordan curves in $\partial N$
which bound a unique least area disk in $\BR^3$. In some sense, this is
generic uniqueness for extremal curves in $\BR^3$.

By using the same setting and similar techniques, Fang-Hua Lin showed
\cite{Li2} that an extremal curve either bounds a unique least area disk,
or two ``extremal'' stable minimal disks. One can get the same result by
using the methods of this paper. Also, by using similar techniques in
an analytical way, Li-Jost proved \cite{Lj} that a $C^{3,\mu}$ Jordan
curve {\em in} $\BH^3$ with total curvature less than $4\pi$ bounds a
unique least area disk.

\subsection{Questions}

As we mentioned in the introduction, there are not many results on the number of solutions of asymptotic Plateau
problem in $\BH^3$. In this paper, we showed a generic uniqueness result for some smooth class of curves. One suspects
whether there is any smooth class of Jordan curves in $\Si$ such that any curve in this class bounds only finitely many
least area planes. Another question in the opposite direction is whether there is any {\it smooth} Jordan curve in
$\Si$ bounding infinitely many least area planes.

If you remove the condition being plane (the topological type of the
surface is disk), Anderson gave examples of Jordan curves bounding
infinitely many area minimizing surfaces in \cite{A2}. If one remove the
embeddedness condition for the curve, the universal cover of a hyperbolic
manifold fibering over a circle, induces a Peano curve in $\Si$ bounding
infinitely many minimal planes, corresponding to the universal cover of
the least area representative of a fiber.

It is still an open question to find a Jordan curve bounding infinitely many least area planes. But it is reasonable to
hope that the bridge principle holds in this case too, and one constructs a rectifiable embedded curve in $\Si$ similar
to the Euclidean case.

\bibliographystyle{gtart}
\bibliography{link}

\end{document}